\documentclass{article}
\usepackage[utf8]{inputenc}

\usepackage[english]{babel}

\usepackage[utf8]{inputenc}
\usepackage[T1]{fontenc}

\usepackage{amssymb}
\usepackage{stmaryrd}
\usepackage{amsmath,amsfonts}
\usepackage[tikz]{bclogo}
\usepackage[top=4.34cm]{geometry}
\usepackage{subcaption}
\usepackage{url}

\usepackage{lmodern}
\usepackage{textcomp}
\usepackage{mathtools, bm}
\usepackage{amssymb, bm}
\usepackage[unq]{unique}
\usepackage{enumerate}
\usepackage{enumitem}
\usepackage{comment}
\usepackage[breaklinks]{hyperref}

\usepackage[numbers]{natbib}  

\usepackage[breaklinks]{hyperref}
\usepackage[numbers]{natbib}

\usepackage{amsthm}

\usepackage{cleveref}

\usepackage{tikz}
\usetikzlibrary {quotes} 

\newtheorem{theorem}{Theorem}[section]
\newtheorem{lemma}[theorem]{Lemma}

\newtheorem{conjecture}[theorem]{Conjecture}

\newtheorem{Problem}[theorem]{Problem}
\newtheorem{Claim}[theorem]{Claim}

\newtheorem{Observation}[theorem]{Observation}

\theoremstyle{definition}

\DeclareMathOperator{\Var}{Var}
\newcommand{\Prob}{\mathbb{P}}
\newcommand{\eps}{\varepsilon}
\newcommand{\st}{\mathbin{\colon}}
\undef{\set}
\DeclarePairedDelimiter{\set}{\lbrace}{\rbrace}
\newcommand{\RadSums}{\mathcal{X}}

\title{Tight lower bounds for anti-concentration of Rademacher sums and Tomaszewski’s counterpart problem}

\date{}

\author{Lawrence Hollom \footnote{\href{mailto:lh569@cam.ac.uk}{lh569@cam.ac.uk}, Department of Pure Mathematics and Mathematical Statistics (DPMMS), University of Cambridge, Wilberforce Road, Cambridge, CB3 0WA, United Kingdom} \and Julien Portier\footnote{\href{mailto:jp899@cam.ac.uk}{jp899@cam.ac.uk}, Department of Pure Mathematics and Mathematical Statistics (DPMMS), University of Cambridge, Wilberforce Road, Cambridge, CB3 0WA, United Kingdom}}

\begin{document}

\maketitle

\begin{abstract} 
    In this paper we prove that $\mathbb{P}(|X| \geq \sqrt{\Var(X)}) \geq 7/32$ for every finite Rademacher sum $X$, confirming a conjecture by Hitczenko and Kwapie{\'n} from 1994, and improving upon results from Burkholder, Oleszkiewicz, and Dvořák and Klein. Moreover we fully determine the function $f(y)= \inf_X \mathbb{P}(|X| \geq y\sqrt{\Var(X)})$ where the $\inf$ is taken over all finite Rademacher sums $X$, confirming a conjecture by Lowther and giving a partial answer to a question by Keller and Klein.
\end{abstract}

\section{Introduction}

In this paper we are interested in tight lower bounds concerning the anti-concentration of finite Rademacher sums. More precisely, if $\mathcal{R}$ is the class of random variables $X$ of the form $X=a_1\eps_1+\dots+a_n\eps_n$ where the $a_i$ are real constants and $\eps_i$ are independent Rademacher variables (that is, picked independently and uniformly at random from $\set{-1,+1}$), we are interested in an exact value of $f(y)=\inf_{X \in \mathcal{R}} \mathbb{P}(|X| \geq y \sqrt{\Var(X)})$ for each $y > 0$.

For $0< y < 1$, Paley-Zygmund inequality already gives 
\begin{align}
\mathbb{P}(|X| \geq y \sqrt{\Var(X)}) = \mathbb{P}(X^2 \geq y^2 \Var(X)) \geq (1-y^2)^2 \frac{\Var(X)^2}{\mathbb{E}[X^4]}
\end{align}
but this is far from giving optimal bounds. For $y>1$, we may have $\mathbb{P}(X \geq y \sqrt{\Var(X)})=0$, for instance for $X=\epsilon_1$. Concerning $y=1$, Burkholder \cite{burkholder1968independent} proved in 1967 that $C =f(1) > 0$.
This was improved to $C \geq e^{-4}/4$ in 1994 by Hitczenko and Kwapie{\'n} \cite{hitczenko1994rademacher}, then to $C \geq 1/10$ in 1996 by Oleszkiewicz \cite{oleszkiewicz1995stein} and, more recently, to $C \geq 3/16$ by Dvořák and Klein \cite{dvorak2022probability}. 
Moreover, Hitczenko and Kwapie{\'n} \cite{hitczenko1994rademacher} conjectured that $C=7/32$, given the tight example $a_1=\dots=a_6$. The problem of determining $C =\inf_{X \in \mathcal{R}} \mathbb{P}(|X| \geq \sqrt{\Var(X)})=2\inf_{X \in \mathcal{R}} \mathbb{P}(X \geq \sqrt{\Var(X)})$ is a natural counterpart of the well known Tomaszewski’s problem \cite{guy1986any} from 1986, which is concerned with determining $\inf_{X \in \mathcal{R}} \mathbb{P}(X \leq \sqrt{\Var(X)})$. This problem attracted significant attention over the years (see for instance \cite{holzman1992product, ben2002robust, boppana2017tomaszewski, boppana2021tomaszewski, dvorak2020improved}) before it was finally settled by Keller and Klein \cite{keller2022proof}, who proved that $\inf_{X \in \mathcal{R}} \mathbb{P}(X \leq \sqrt{\Var(X)})=1/2$.  \\

Inspired by a computer simulation, Lowther \cite{Lowther} conjectured that the function  $f(y)=\inf_{X} \mathbb{P}(|X| \geq y \sqrt{\Var(X)})$ is piecewise constant on $\mathbb{R}_{\geq 0}$, taking seven different values, with discontinuities at 0, $1/\sqrt{7}$, $1/\sqrt{5}$, $1/\sqrt{3}$, $2/\sqrt{6}$ and $1$. In this paper we prove Lowther's conjecture.

\begin{theorem}
\label{thm:ExplicitInfFunction}
    Let the function $f: \mathbb{R}_{\geq 0} \rightarrow \mathbb{R}$ be defined for every $y \geq 0$ by $f(y)=\inf_{X \in \mathcal{R}} \mathbb{P}(|X| \geq y \sqrt{\Var(X)})$, where the inf is taken over the class $\mathcal{R}$ of all finite Rademacher sums. Then the function $f$ is piecewise constant. More precisely,
    
    \begin{itemize}
        \item $f(0)=1$,
        \item $f(y)=1/2$ for $y \in (0,1/\sqrt{7}]$,
        \item $f(y)=29/64$ for $y \in (1/\sqrt{7},1/\sqrt{5}]$,
        \item $f(y)=3/8$ for $y \in (1/\sqrt{5},1/\sqrt{3}]$,
        \item $f(y)=1/4$ for $y \in (1/\sqrt{3},2/\sqrt{6}]$,
        \item $f(y)=7/32$ for $y \in (2/\sqrt{6},1]$,
        \item $f(y)=0$ for $y \in (1, \infty)$.
    \end{itemize}
\end{theorem}

In particular, \Cref{thm:ExplicitInfFunction} shows that $\mathbb{P}(|X| \geq \sqrt{\Var(X)}) \geq 7/32$, which proves the conjecture of Hitczenko and Kwapie{\'n} \cite{hitczenko1994rademacher}.

\section{Proof outline}

Our proof follows similar lines to that of Dvořák and Klein \cite{dvorak2022probability} of the weaker bound $\mathbb{P}(|X| \geq \sqrt{\Var(X)}) \geq 3/16$, and many of our techniques are closely related to theirs. 
We now introduce some notation, and then proceed to outline our proof, indicating where it differs from the method of Dvořák and Klein.

\subsection{Notation}
\label{subsec:notation}

Throughout this paper, $X$ will be a sum of Rademacher random variables, i.e. $X=\sum_{i=1}^n a_i\eps_i$, where $a_i$ are real constants satisfying $a_1\geq a_2\geq\dots\geq a_n >0$, and each $\eps_i$ is a random variable picked uniformly from $\set{-1,+1}$.
It will be convenient to assume that $\Var(X)=1$, or equivalently that $a_1^2+\dots+a_n^2=1$, so we let $\RadSums$ be the set of all such Rademacher sums.

We will often end up in a situation where we know that, for some $k<n$, each of the $(a_i)_{1\leq i\leq k}$ lie within some small interval.
In such a case, we will write $a_i=c_i+\delta_i$, where each $c_i$ will be a known constant, and we will have bounds of the form $|\delta_i|\leq \delta$ for some constant $\delta$ and each $1\leq i \leq k$.

Furthermore, in cases like the above, we will in practise always have that $c_1^2+\dots+c_k^2=1$, and so we will know that $a_{k+1},\dotsc,a_n$ will be small.
This partition into large and small $a_i$ will be represented by writing $X=Y_k+Z_k$, where $Y_k=a_1\eps_1+\dots+a_k\eps_k$ and $Z_k=a_{k+1}\eps_{k+1}+\dots+a_n\eps_n$.

\subsection{Algorithm overview}
\label{subsec:algorithm-overview}

A key tool in our proof is the use of a computerised bound of $\tilde{G}(a,x)=\inf_{X \in \RadSums_a} \mathbb{P}(X \geq x)$ where $\RadSums_a$ is the class of Rademacher sums of variance 1 satisfying $a_1 \leq a$.

We lower bound this function $\tilde{G}(a,x)$ by a function $D(a,x)$, which can itself be lower bounded in a manner which can be efficiently computed.
We explain this bounding process further in \Cref{subsec:computer-simulation}, but note here that our bounding function $D(a,x)$ is identical to that of Dvořák and Klein \cite{dvorak2022probability}.

What is new in our method is how we use this lower bound $D(a,x)$ to automatically sift through cases.
The key tool in this automated process is the elimination of variables, by which we may remove the $k$ largest weights from a Rademacher sum, to be left with only the tail to consider.
See \Cref{obs:elimination} for a precise statement of this elimination of variables.

With these tools in hand, we may now outline our algorithm. 
Whereas Dvořák and Klein proceeded via taking a mesh over possible values for $a_1,a_2,a_3$, we proceed in an iterative manner.
This allows us to, for instance, skip considering all possible values for $a_2,a_3,\cdots,a_k$ when we can already conclude our desired result from the value of $a_1$ alone.

Our algorithm first considers which values $a_1$ might take, splitting the interval $[0,1]$ into some number $d$ of smaller intervals of width $1/d$.
For each of these ranges of values of $a_1$, we consider eliminating $a_1$ from the sum. 
In some cases, this allows us to apply our bounding function $D(a,x)$ and immediately conclude the desired result, which will be of the form $\Prob(X\geq s)\geq p$ for some $s$ and $p$.
We discard all such intervals.

For each interval which is not discarded, we repeat the process for $a_2$, considering which values it may take, eliminating $a_1$ and $a_2$, and applying our bounding function.
We continue this process for $a_3,\dotsc,a_k$, for some fixed $k$.
When the algorithm terminates, we will be left with a list of $k$-tuples of intervals of values for $a_1,\dotsc,a_k$ which have not been dealt with.
By choosing suitably large $d$ and $k$, and applying the above algorithm iteratively (see \Cref{subsubsec:ImprovedAlgoBoundingai} for a more thorough explanation and \Cref{appendix:simulation} for the precise details), we can reduce our problem to a small enough number of regions that they may be checked by hand; these are the neighbourhoods of the `hard cases'.

\subsection{Neighbourhoods of hard cases}
\label{subsec:nbhds-of-hard-cases}

The term `hard case' refers to a sequence of $a_i$'s for which we have $\Prob(X\geq s)\geq p$, whereas $\Prob(X > s) < p$.
These cases cannot be dealt with by our algorithm, as our bounding function $D(a,x)$ is in fact also a lower bound for $\Prob(X > x)$.
Therefore, we cannot hope for these hard cases and their neighbourhoods to be dealt with by the algorithm, and so they must be dealt with by another method.

Our second improvement is a general and simple method to deal with these situations, provided that the neighbourhood is sufficiently small.
We do this by means of two key results, \Cref{lem:LargeVarianceLeftNew} and \Cref{lem:small-sum}.
If we let $c_1,\dotsc,c_k$ be the hard case itself, so $a_i$ is close to $c_i$ for $1\leq i \leq k$, then these lemmas intuitively deal with the following two cases.
We write $X=Y_k+Z_k$, where $Y_k$ captures the first $k$ terms of $X$, and $Z_k$ is the tail.
If $\Var(Z_k)$ is on the order of $\max_i|a_i-c_i|$, then whenever $Y_k\approx s$ we can bound $\Prob(X\geq s)$ away from zero, which will allow us to conclude.
Conversely, if $\Var(Z_k)$ is small, then we know the values of the $a_i$ cannot be too small, and we will be able to show that in many cases where $Y_k\approx s$, we in fact have $Y_k\geq s$, which will also allow us to conclude.

In the next section, we give full details of the tools we have mentioned in this proof outline.

\section{Tools}

\subsection{Combinatorial and Probabilistic Tools}

The first tool we present allows us to eliminate some of the $a_i$ from a Rademacher sum.

\begin{Observation}
\label{obs:elimination}
    If $X=Y_k+Z_k$ is a Rademacher sum, where $Y_k=\sum_{i=1}^ka_i\eps_i$ and $Z_k=\sum_{i=k+1}^n a_i\eps_i$, and $s$ is some real constant, then we have the following.
    \begin{align*}
        \Prob(X\geq s)=2^{-k}\sum_{\zeta\in\set{-1,+1}^k}\Prob\Bigl[Z_k\geq s+\sum_{i=1}^ka_i\zeta_i\Bigr]
    \end{align*}
\end{Observation}

The proof of \Cref{obs:elimination} is straightforward, and a complete discussion can be found in \cite[Lemma 2.1]{keller2022proof}.

We will use the following observation from \cite{dvorak2022probability}, which is a consequence of a well known result of Erd{\"o}s \cite{erdos1945lemma}.

\begin{Observation}
\label{obs:observation1DK}
Let $b_1 \geq b_2 \geq \dots \geq b_t >0$ be such that $b_{t-k+1} + \dots + b_t \geq \alpha$ for some $\alpha >0$ and $t \geq k > 0$. Then for any $x$ and any $b_{t+1}, \dots, b_t$ we have 

$$\mathbb{P}[ \sum_{i=1}^{s} b_i\eps_i \in (x-\alpha, x+\alpha)] \leq f(k,t)/2^t$$
where $f(k,t)$ denotes the sum of the $k$ largest binomial coefficients of the form $\binom{t}{i}$ for some $i$, $0 \leq i \leq t$.
\end{Observation}

We will also use the following result \cite[Theorem 1.3]{dvorak2022probability}.

\begin{theorem}
\label{thm:WeakLowther0.35}
    For every Rademacher sum $X$, we have $\mathbb{P}(X \geq 0.35\sqrt{\Var(X)}) \geq 1/4$.
\end{theorem}

We now present our technical lemmas that we will use many times throughout our proof.
The most important of these are \Cref{lem:LargeVarianceLeftNew} and \Cref{lem:small-sum}, which deal with the two possibilities discussed in \Cref{subsec:nbhds-of-hard-cases}, depending on how much variance is in the tail of a Rademacher sum. We now give precise statements and proofs of these lemmas.

The first lemma deals with the case when $\Var(Z)$ is large.

\begin{lemma}
\label{lem:LargeVarianceLeftNew}
    Let $s, \gamma, d$ be real numbers with $0 < \gamma < d/2$, and $X=Y+Z$, where $Y$ and $Z$ are independent random variables, and $Z$ is a Rademacher sum with $\Var(Z)\geq 8.17\gamma^2$.
    Define, for integer $r$, the probability $p_r = \Prob(Y\in [s + rd -\gamma, s + rd + \gamma])$, and let $t=\sum_{r=2}^\infty p_r$.
    Then, if $p_{-1}\geq t$ and $p_0\geq p_1$, then $\Prob(X\geq s)\geq p_0/4 + 3p_1/4 + t$.
\end{lemma}
\begin{proof}
    Let $q_1=\Prob(Z\in [\gamma, d+\gamma))$ and $q_2=\Prob(Z\geq d+\gamma)$. Noting that $8.17\geq 1/0.35^2$ and using \Cref{thm:WeakLowther0.35}, we have
    \begin{align*}
        q_1+q_2 = \mathbb{P}(Z \geq \gamma) \geq \mathbb{P}(Z \geq 0.35\sqrt{\Var(Z)}) \geq 1/4.
    \end{align*}

    Then by applying the assumptions in the statement we see that
    \begin{align*}
        \Prob(X\geq s) &\geq \sum_{r=-1}^\infty p_r\Prob(Z\geq \gamma - rd) \\
        &\geq p_{-1}\Prob(Z\geq \gamma+d) + p_0\Prob(Z\geq \gamma) + p_1\Prob(Z\geq -\gamma) + t\Prob(Z\geq -d-\gamma) \\
        &\geq p_{-1}q_2 + p_0(q_1+q_2) + p_1(1 - q_1 - q_2) + t(1 - q_2) \\
        &= (q_1+q_2)(p_0 - p_1) + q_2(p_{-1} - t) + p_1 + t \\
        &\geq p_0/4 + 3p_1/4 + t.
    \end{align*}
    As required.
\end{proof}

We would like to remark that we will eventually prove \Cref{thm:sqrt-7-thm}, which is a strengthening of \Cref{thm:WeakLowther0.35}. Hence, by using \Cref{thm:sqrt-7-thm} instead of \Cref{thm:WeakLowther0.35} in the proof of \Cref{lem:LargeVarianceLeftNew}, one can relax the lower bound on $\Var(Z)$ to $\Var(Z) \geq 7\gamma^2$. However, since we use \Cref{lem:LargeVarianceLeftNew} to prove \Cref{thm:sqrt-7-thm}, we will stick to the above version of \Cref{lem:LargeVarianceLeftNew} in this paper. \\

We now state and prove the following `pairing' lemma.

\begin{lemma}
\label{lem:pairing}
    Let $s, a_1,\dotsc, a_k$ be positive reals, and let $A,B\subseteq\set{-1,+1}^k$ be disjoint sets of sequences of signs.
    Assume that we have an injection $\phi\st A\to B$ satisfying
    \begin{align}
    \label{eq:pairing-assumption}
        \sum_{i=1}^k\zeta_i a_i +\sum_{i=1}^k\phi(\zeta)_ia_i \geq 2s
    \end{align}
    for all $\zeta\in A$. Then for any symmetric random variable $Z$, we have that
    \begin{align*}
        \sum_{\eps\in A\cup B}\Prob\Bigl[\sum_{i=1}^k\eps_i a_i + Z\geq s\Bigr]\geq |A|.
    \end{align*}
\end{lemma}
\begin{proof}
    For any $\zeta\in\set{-1,+1}^k$, let $g(\zeta)=\sum_{i=1}^k\zeta_i a_i$. We can then compute:
    \begin{align*}
        \sum_{\zeta\in A\cup B}\Prob(g(\zeta) + Z\geq s) &\geq \sum_{\zeta\in A}\bigl(\Prob(g(\zeta) + Z\geq s) + \Prob(g(\phi(\zeta)) + Z\geq s)\bigr) \\
        &= \sum_{\zeta\in A}\bigl(\Prob(Z\geq s - g(\zeta)) + 1 - \Prob(Z\geq g(\phi(\zeta)) - s)\bigr)\\
        &\geq |A|.
    \end{align*}
    Where the final inequality follows by application of \eqref{eq:pairing-assumption}.
\end{proof}

Using \Cref{lem:pairing}, we prove our second technical lemma, which intuitively will be applied in the cases where $\Var(Z)$ is small.

\begin{lemma}
\label{lem:small-sum}
    Assume that we have numbers $c_i$ and $\delta > 0$ such that $a_i\in[c_i-\delta,c_i+\delta]$ for each $i$.
    Assume further that $X\in\RadSums$, and we have constants $s$ and $p$ such that $\Prob(X\geq s) < p$.
    Take some $I\subseteq [k]$, and fix $\lambda_i\in\set{-1,+1}$ for each $i\in I$, so that $\sum_{i\in I} \lambda_i c_i = s$.
    Then define the following three sets.
    \begin{align*}
        S &= \set{\zeta\in\set{-1.+1}^k\st \sum_{i=1}^kc_i\zeta_i = s} \\
        R &= \set{\zeta\in S\st \zeta_i = \lambda_i \text{ for each }i \in I} \\
        T &= \set{\zeta\in\set{-1.+1}^k\st \sum_{i=1}^kc_i\zeta_i > s}.
    \end{align*}
    Let $d$ be such that $\min\set{\sum_{i=1}^k c_i\zeta_i\st \zeta\in T} = s+d$.
    Assume that we have that $\delta \leq d/k$ and 
    \begin{align}
    \label{eq:small-sum-condition}
        \min(|S|-|T|,|R|)/2 + |T| \geq 2^kp.
    \end{align} 
    Then we have that $\sum_{i\in I}\lambda_i a_i < s$.
\end{lemma}
\begin{proof}
    Assume for contradiction that $\sum_{i\in I} \lambda_i a_i \geq s$.
    We can write $R$ as a union of pairs of the form $\zeta, \eta$, where $\eta$ is simply $\zeta$ flipped on $[k]\setminus I$.
    Then if we define $g(\eps)=\sum_{i=1}^k a_i\eps_i$, we know that $g(\zeta)+g(\eta) = 2\sum_{i\in I} \lambda_i a_i \geq 2s$, and thus by \Cref{lem:pairing} we may conclude that 
    \begin{align*}
        \sum_{\eps\in R}\Prob\Bigl[Z\geq s - \sum_{i=1}^k a_i\eps_i\Bigr] \geq |R| / 2.
    \end{align*}
    Now take a set $U\subseteq R$ satisfying both $|U| = \min(|S|-|T|,|R|)$ and
    \begin{align}
    \label{eq:evaluate-U-prob-sums}
        \sum_{\eps\in U}\Prob\Bigl[Z\geq s - \sum_{i=1}^k a_i\eps_i\Bigr] \geq |U| / 2.
    \end{align}
    Note that this may be achieved by deleting from $R$ the elements for which the above summand is minimal.

    Now, noting that $|S| - |U|\geq |T|$, we may pick an arbitrary injection $\phi:T\rightarrow S\setminus U$.
    For $\zeta\in T$, we deduce from $|\delta_i|\leq \delta\leq d/k$ that $|s - \sum_{i=1}^k a_i\zeta_i| \leq |s - \sum_{i=1}^k a_i \phi(\zeta)_i|$.
    Therefore by \Cref{lem:pairing} again we have that
    \begin{align}
    \label{eq:evaluate-paired-T-prob-sums}
        \sum_{\eps\in T}\bigl[\Prob(Z\geq s - \sum_{i=1}^k a_i\eps_i) + \Prob(Z\geq s - \sum_{i=1}^k a_i\phi(\eps)_i)\bigr] \geq |T|.
    \end{align}

    Thus by eliminating $a_1,\dotsc,a_k$ from $X$, we find by \Cref{obs:elimination} that 
    \begin{align*}
        p &> \Prob(X\geq s) \\
        &\geq 2^{-k}\biggl(\sum_{\eps\in U}\Prob\Bigl[Z\geq s - \sum_{i=1}^ka_i\eps_i\Bigr] + \sum_{\eps\in T}\Prob\Bigl[Z\geq s - \sum_{i=1}^ka_i\eps_i\Bigr] + \sum_{\eps\in S\setminus U}\Prob\Bigl[Z\geq s - \sum_{i=1}^ka_i\eps_i\Bigr]\biggr) \\
        &\geq 2^{-k}(|U|/2 + |T|) \\
        &\geq p.
    \end{align*}
    Where the penultimate inequality comes from applying inequalities \eqref{eq:evaluate-U-prob-sums} and \eqref{eq:evaluate-paired-T-prob-sums}. This is a contradiction, and so we find that $\sum_{i\in I}\lambda_i a_i < s$, as required.
\end{proof}

In practise, the key condition which we will need to check when applying the above lemma is \eqref{eq:small-sum-condition}.
Note that in every case we will actually apply the lemma, we will take the $c_i$ all to be integer multiples of $c_k$, assuming that the $c_i$ are non-increasing. 
In this case then the value of $d$ in the above proof can be bounded below by $2c_k$, so it suffices to check that $\delta \leq 2c_k/k$. \\

The following simple inequality will also be of use.

\begin{lemma}
\label{lem:IneqSumSquares}
Let $k \in \mathbb{N}$ and $a_1, \dots, a_k$, $b$ and $\eps$ be some reals such that $\sum_{i=1}^k a_i \leq b$ and $b/k-\eps \leq a_1, \dots, a_k \leq b/k+\eps$. Then $\sum_{i=1}^k a_i^2 \leq b^2/k+k\eps^2$.
\end{lemma}

\begin{proof}
    Let $a_i=b/k+\delta_i$. Note that by the assumptions we have $\sum_{i=1}^k \delta_i \leq 0$ and $|\delta_i| \leq \eps$ for all $i$. By simple rearranging we have that
    \begin{align*}
        \sum_{i=1}^k a_i^2 &=\sum_{i=1}^k (b/k+\delta_i)^2 \\
        &= \sum_{i=1}^k b^2/k^2+2b\delta_i/k+\delta_i^2 \\
        &\leq b^2/k+ k\eps^2.
    \end{align*}
\end{proof}

\subsection{Computer-aided bounds}
\label{subsec:computer-simulation}

The goal of this subsection is to explain the computer-aided bounds. As mentioned in the proof outline, we first find lower bounds $D(a,x)$ of $\inf_{X \in \RadSums_a} \mathbb{P}(X \geq x)$ using a computer-aided computation exactly as Dvořák and Klein \cite{dvorak2022probability}. As it is an important part of our proof and for sake of completeness, we re-explain the details in \ref{subsubsec:DynamicProgrammingBound}. We then explain the details of our improved algorithm giving bounds on the $a_i$ in \ref{subsubsec:ImprovedAlgoBoundingai}.

\subsubsection{Dynamic programming bound}
\label{subsubsec:DynamicProgrammingBound}

We will iteratively compute the functions
\begin{align*}
    D_i : (0,1) \times \mathbb{R} \rightarrow \mathbb{R}
\end{align*}
that will satisfy $\tilde{G}(a, x) \geq D_i(a,x)$ for every $i$, $a$ and $x$. We first initialise $D_0(a,x)$ to the lower bound given by Prawitz's inequality \cite{prawitz1972limits}. More precisely, in the case of Rademacher sums, the following formula was made explicit in \cite{keller2022proof} and used in \cite{dvorak2022probability}.

\begin{lemma}
\label{lem:Prawitz}
    For all $q \in [0,1]$ and  $T >0$, we have $\tilde{G}(a_1, x) \geq F(a_1, x, T, q)$, where

    \begin{equation}
        F(a,x,T,q)=1/2- \int_{0}^{q} |k(u,x,T)|g(Tu,a)du - \int_{q}^{1} |k(u,x,T)|h(Tu,a)du - \int_{0}^{q} k(u,x,T)\exp(-(Tu)^2/2)du
    \end{equation}

    \begin{equation}
        k(u,x,T)=\frac{(1-u)\sin(\pi u+Tux)}{\sin(\pi u)}+\frac{\sin(Tux)}{\pi}
    \end{equation}

    \begin{equation}
  g(v,a) =
    \begin{cases}
      \exp(-v^2/2)-\cos(av)^{1/a^2} & \text{if } av \leq \pi/2 \\
      \exp(-v^2/2)+1 & \text{otherwise} \\
    \end{cases}   
\end{equation}

\begin{equation}
  h(v,a) =
    \begin{cases}
      \exp(-v^2/2) & \text{if } av \leq \theta \\
      (-\cos(av))^{1/a^2} & \text{if } \theta \leq av \leq \pi \\
      1 & \text{otherwise} \\
    \end{cases}   
\end{equation}

and $\theta = 1.778 \pm 10^{-4}$ is the unique solution of $\exp(-\theta^2/2) = - \cos(\theta)$ in the interval $[0, \pi]$.
\end{lemma}

In regard to the induction, note that by elimination of $a_1$ we have the following.
\begin{equation}
\label{eq:BoundEliminationa1}
    \tilde{G}(a^*,x)= \frac{1}{2} \inf_{a \in (0,a^*]} \Bigl\{ \tilde{G}\bigl(\frac{a}{\sqrt{1-a^2}}, \frac{x-a}{\sqrt{1-a^2}}\bigr) + \tilde{G}\bigl(\frac{a}{\sqrt{1-a^2}}, \frac{x+a}{\sqrt{1-a^2}}\bigr) \Bigr\}.
\end{equation}

Hence, we compute recursively the $D_i$ using the following formula
\begin{equation}
\label{eq:RecursionD}
    D_{i+1}(a^*,x)= \max \Bigl[D_{i}(a^*,x), \frac{1}{2} \inf_{a \in (0,a^*]} \Bigl\{ D_i\bigl(\frac{a}{\sqrt{1-a^2}}, \frac{x-a}{\sqrt{1-a^2}}\bigr) + D_i\bigl(\frac{a}{\sqrt{1-a^2}}, \frac{x+a}{\sqrt{1-a^2}}\bigr) \Bigr\}\Bigr].
\end{equation}

It follows by induction and \eqref{eq:BoundEliminationa1} that $\tilde{G}(a^*,x) \geq D_i(a^*,x)$ for every $i$. In practice, we compute $D_i(a^*, x)$ for $a^* \in [0, 1]$ and $x \in [-3, 3]$ with granularity of $\beta = 1/2000$ (with $a^*$ starting from $0$ and $x$ starting from $-3$). We replace \eqref{eq:RecursionD} with a discrete variant where arguments of $D_{i+1}$ are rounded up to a multiple of $\beta$, hence underestimating $D_{i+1}$, and we also apply this rounding-up to the $\frac{a}{\sqrt{1-a^2}}$ and the $\frac{x \pm a}{\sqrt{1-a^2}}$ arguments. Hence, we only consider a finite set of $a \in [0, a^*]$ in the infimum at \eqref{eq:RecursionD} to compute $D_{i+1}(a^*, x)$, and this results in a dynamic-programming method for computing the $D_i(a^*, x)$'s. Using this method, we compute $D_I(a,x)$ for $I=1000$, we set $D(a,x)=D_I(a,x)$, and $D(a,x)$ is the lower bound we will use for $\tilde{G}(a,x)$ in the rest of this paper. 

\subsubsection{Iterative bounding algorithm}
\label{subsubsec:ImprovedAlgoBoundingai}

As explained in the proof outline, for fixed $s$, $p$ and $k$ our algorithm finds intervals of $[0,1]^k$ in which a candidate counterexample to $\mathbb{P}(X \geq s)<p$ must lie.
For instance, in \Cref{sec:tomaszewski-counterpart} we will be mostly interested in the case $s=1$ and $p=7/64$. We now explain this algorithm in greater detail. 

We first fix some positive integer $d$ and divide up the interval $[0,1]$ into $d$ intervals of width $1/d$. We then test, using the function $D$, which of these intervals contain a value for $a_1$ which could be the start of a sequence for which $\mathbb{P}(X \geq s)<p$. 
More precisely, for each interval $[a_{1,-},a_{1,+}]$, let $\sigma=\sqrt{1-a_{1,+}^2}$ be the minimum possible standard deviation of the rest of the Rademacher sum.
We then eliminate $a_1$ by \Cref{obs:elimination} and check whether
\begin{equation}
\label{eq:algorithm-test-one-variable}
    \frac{1}{2}D\Bigl(\frac{a_{1,+}}{\sigma}, \frac{s-a_{1,-}}{\sigma}\Bigr)+\frac{1}{2}D\Bigl(\frac{a_{1,+}}{\sigma}, \frac{s+a_{1,+}}{\sigma}\Bigr) \geq p.
\end{equation}

If it is, then we discard this interval $[a_{1,-},a_{1,+}]$. 
If it is not, then we reiterate the same process for $a_2$: split $[0,1]$ into intervals and test by elimination of $a_1$ and $a_2$ which pairs of intervals $[a_{1,-},a_{1,+}] \times [a_{2,-},a_{2,+}]$ we can discard.
This process is continued until some fixed $a_k$, keeping track of the interval of values of $a_i$ which cannot be ruled out for each $i$, at which point the process halts.

To generalise equation \eqref{eq:algorithm-test-one-variable} to more variables, assume we are considering $a_1,\dotsc,a_r$, and in particular whether they may take values in $[a_{1,-},a_{1,+}]\times\dots\times[a_{r,-},a_{r,+}]$.
Set 
\begin{align*}
    \sigma_-=\sqrt{1-a_{1,+}^2-\dots-a_{r,+}^2}\text{ and }\sigma_+=\sqrt{1-a_{1,-}^2-\dots-a_{r,-}^2},
\end{align*}
and assume for now that these are both real.
We now define
\begin{equation}
    h(\zeta)=\max_{z_i\in[a_{i,-},a_{i,+}]}\sum_{i=1}^k \zeta_i z_i.
\end{equation}

Then, noting that the above maximum will take $z_i=a_{i,-}$ if $\zeta_i=-1$ and $z_i=a_{i,+}$ if $\zeta_i=+1$, we need to check whether 

\begin{equation}
\label{eq:algorithm-test-many-variables}
    2^{-r}\sum_{\zeta\in\set{-1,+1}^r}D\Bigl(\frac{a_{r,+}}{\sigma_-}, \max_{\sigma\in\set{\sigma_-,\sigma_+}}\frac{s+h(\zeta)}{\sigma}\Bigr)\geq p.
\end{equation}

We may note that the maximum in equation \eqref{eq:algorithm-test-many-variables} will take $\sigma=\sigma_-$ if $s+h(\zeta)\geq 0$ and $\sigma=\sigma_+$ otherwise.

However, we run into an issue if $\sum_{i=1}^ra_{i,-}^2 \leq 1 < \sum_{i=1}^r a_{i,+}^2$, as then $\sigma_-$ is not real. (Note that if $\sum_{i=1}^ra_{i,-}^2>1$ then we may immediately discard the sequence as it contradicts $\Var(X)=1$.)
In this case, all terms with $s+h(\zeta)\geq 0$ contribute nothing to the sum, as we cannot bound the remaining variance away from zero.
Thus the condition that we check in this case is

\begin{equation}
\label{eq:algorithm-test-no-variance-left}
    2^{-r}\sum_{\substack{\zeta\in\set{-1,+1}^r \\ x+h(\zeta) < 0}}D\Bigl(1, \frac{s+h(\zeta)}{\sigma_+}\Bigr)\geq p.
\end{equation}

The above algorithm produces a collection $B$ of $k$-tuples of intervals such that we have bounds of the form $(a_1, \dots, a_k) \in \bigcup_{b\in B} [b_1^{-}, b_1^{+}] \times \dots \times [b_k^{-}, b_k^{+}]$ for any candidate $X=\sum_{i=1}^{n} a_i\eps_i$ satisfying $\mathbb{P}(X \geq s) < p$. 
Note that these bounds are rigorous bounds, and not merely approximations, as all the functions of the $a_i$'s we are dealing with in this computation are monotone in the domains we consider them over. 
For this reason, in the computation, we always bound $a_i$ by the upper bound or lower bound of their interval, depending on whether the function is increasing or decreasing.

The final step of our algorithm is to output the bounds on each $a_i$ which it can prove.
These can then be fed back into the start of the algorithm; if, for example, we have shown on one run of the algorithm that $0.24\leq a_5\leq 0.26$ then we never need to check values outside this range.
This allows us to re-run the algorithm with a larger value of $d$, which in turn allows us to produce more precise bounds without an excessively long computation.

The precise details of the operation of the code can be found in the code itself \cite{simulation}. Moreover, a list of all results proved by the computer program can be found in \Cref{appendix:simulation}.
\section{Tomaszewski's counterpart}
\label{sec:tomaszewski-counterpart}

In this section we prove the following result, which can be seen as a tight bound for the counterpart Tomaszewski problem.

\begin{theorem}
\label{thm:main-thm}
    For any Rademacher sum $X=\sum_{i=1}^n a_i\eps_i\in\RadSums$, we have
    $$\Prob(|X|\geq 1)\geq 7/32.$$
\end{theorem}

We first point out the two following results, which follow by a direct application of \Cref{obs:observation1DK}.

\begin{lemma}
\label{lem:a1+a2<1}
    If $\Prob(X\geq 1) < 7/64$ then $a_1+a_2 < 1$.
\end{lemma}

\begin{lemma}
\label{lem:sum_a3a4a5}
    If $\Prob(X\geq 1) < 7/64$ then $a_3+a_4+a_5 < 1$.
\end{lemma}

Those two lemmas will speed up the automated computation, but also be used in our proofs to take care of the neighbourhood of the hard cases. \\

By using our automated computation, we can reduce to the following eight cases.

\begin{lemma}
Suppose $X=\sum_{i=1}^n a_i\epsilon_i \in \mathcal{X}$  satisfies $\mathbb{P}(X \geq 1) < 7/64$. Then we must be in one of the following cases.

\begin{enumerate}[label=(\Alph*)]
    \item $a_1\geq 1-\delta$, where $\delta = 0.04$.
    \item $a_1,\dotsc,a_4\in (1/2 - \delta,1/2+\delta)$, where $\delta = 0.005$.
    \item $a_1,\dotsc,a_9\in (1/3 - \delta, 1/3 + \delta)$, where $\delta = 0.07$.
    \item $a_1,\dotsc,a_{16}\in(1/4 - \delta, 1/4 + \delta)$, where $\delta = 0.03$.
    \item $a_1\in (2/3-\delta, 2/3+\delta)$ and $a_2,\dotsc,a_6\in (1/3 - \delta, 1/3 + \delta)$, where $\delta = 0.002$.
    \item $a_1\in (1/2-\delta, 1/2+\delta)$ and $a_2,\dotsc,a_{13}\in (1/4 - \delta, 1/4 + \delta)$, where $\delta = 0.03$.
    \item $a_1,a_2\in (1/2-\delta, 1/2+\delta)$ and $a_3,\dotsc,a_{10}\in (1/4 - \delta, 1/4 + \delta)$, where $\delta = 0.05$.
    \item $a_1,a_2,a_3\in (1/2-\delta, 1/2+\delta)$ and $a_4,\dotsc,a_7\in (1/4 - \delta, 1/4 + \delta)$, where $\delta = 0.005$.
\end{enumerate}
\end{lemma}

\begin{proof}
    By automated computation. See \Cref{appendix:simulation} item \ref{sim:0} and \ref{sim:0A} - \ref{sim:0L}.
\end{proof}

We shall consider these cases in turn. Some of them were already solved by Dvořák and Klein for some very small $\delta$, namely cases B, E, G and H for $\delta=10^{-9}$ (Proposition 6.2 in \cite{dvorak2022probability}) but the computer-aided bounds on $\delta$ are far from this order. Hence, using a different method than Dvořák and Klein, we will treat those cases for a significantly larger value of $\delta$.
Moreover, in cases C, D, F, and G, we use the value of $\delta$ shown above to prove some instance of \Cref{lem:small-sum}, and then feed the resulting inequality back into the automated computation, producing tighter bounds in a reasonable time frame. 
This allows us to compute a significantly smaller value of $\delta$, which is then used in later parts of the proof. \\

We now consider the above cases in turn, supposing for contradiction in each case that we have a counterexample to \Cref{thm:main-thm}.
Throughout the proof, we will use the notation discussed in \Cref{subsec:notation} without further comment.

\subsection{Case A}
\label{subsec:case-A}
This follows from \cite[Proposition 6.1]{dvorak2022probability}.
The proof given there allows for any $\delta \leq 1/15$, and our value is easily within this bound.

\subsection{Case B}
\label{subsec:case-B}

Suppose for contradiction that there exists a counterexample to \Cref{thm:main-thm} with $a_i=1/2+\delta_i$ and $|\delta_i| \leq \delta \leq 0.005$ for $1 \leq i \leq 4$. Let $\Delta=\max \{|\delta_1|, \dots, |\delta_4|\}$. By \Cref{lem:a1+a2<1}, we have $a_1+a_2 < 1$ and consequently $\Delta=-\delta_4$. By \Cref{lem:IneqSumSquares}, we have $a_1^2+a_2^2 \leq 1/2+2\Delta^2$. Hence
\begin{align*}
    a_1^2+a_2^2+a_3^2+a_4^2 &\leq 1/2+2\Delta^2+1/4+(1/2-\Delta)^2 \\
    &= 1-\Delta+3\Delta^2.
\end{align*}

Letting $Y_4=\sum_{i=1}^4 a_i\epsilon_i$ and $Z_4=\sum_{i=5}^n a_i\epsilon_i$, it follows that $\Var(Z_4)=1-a_1^2-a_2^2-a_3^2-a_4^2 \geq \Delta-3\Delta^2$. 
Let $s=1$, $\gamma=4\Delta$, and $d=1$. Using the notations of \Cref{lem:LargeVarianceLeftNew} with $Y=Y_4$, we have
\begin{align*}
    p_{-1}=\Prob(Y_4\in[-4\Delta,4\Delta])&=\binom{4}{2}/2^4=6/16, \\
    p_0=\Prob(Y_4\in[1-4\Delta,1+4\Delta])&=\binom{4}{3}/2^4=4/16, \\
    p_1=\Prob(Y_4\in[2-4\Delta,2+4\Delta])&=\binom{4}{4}/2^4=1/16, \\
    t=\sum_{r=2}^\infty p_r&=0.
\end{align*}
Note that $p_{-1}\geq t$ and $p_0\geq p_1$.
Moreover, as $\Delta\leq 0.005$ (by \Cref{appendix:simulation} \cref{sim:0B}), $8.17\gamma^2=130.72\Delta^2$, and $\Delta-3\Delta^2 \geq 130.72\Delta^2$ holds for all $\Delta\leq 0.0074$, we may apply \Cref{lem:LargeVarianceLeftNew} and conclude that
\begin{align*}
    \Prob(X \geq 1) &\geq p_0/4+3p_1/4+t
    = 7/64,
\end{align*}
as desired.

\subsection{Case C}
\label{subsec:case-C}

Suppose for contradiction that there exists a counterexample to \Cref{thm:main-thm} with $a_i=1/3+\delta_i$ and $|\delta_i| \leq \delta \leq 0.07$ for $1 \leq i \leq 9$.

\begin{Claim}
\label{claim:case-c-claim-1}
    We have $a_1+a_2+a_3 < 1$, and $\delta < 0.0009$.
\end{Claim}

\begin{proof}
    It suffices to verify that the assumptions of \Cref{lem:small-sum} do hold for $c_i=1/3$ for $1 \leq i \leq 9$, $I=\set{1,2,3}$, $\lambda_i = 1$ for all $i\in I$, and $s=1$. 
    
    Indeed, $\delta\leq 0.07 < 2(1/3)/9$, as required.
    Here $|S| = \binom{9}{3} = 84$, $|R| = \binom{6}{3} = 20$, and $|T| = \binom{9}{2} + \binom{9}{1} + \binom{9}{0} = 46$.
    Thus $\min(|S|-|T|,|R|)/2 + |T| = 56 \geq 2^9(7/64)$,
    and so by \Cref{lem:small-sum} we have $a_1+a_2+a_3 < 1$.
    We may now add this hypothesis into the automated computation (see \Cref{appendix:simulation} \cref{sim:0C}), from which we get that $\delta < 0.0009$, as required.
\end{proof}

Let $\Delta=\max \{|\delta_1|, \dots, |\delta_9|\}$. Note that $\Delta=\max\{|\delta_1|, -\delta_9 \}$. By \Cref{lem:IneqSumSquares}, we have
\begin{align*}
a_1^2+a_2^2+a_3^2 &\leq 1/3+3\Delta^2
\end{align*}

Moreover, $a_3 \leq 1/3-|\delta_1|/2$. Hence
\begin{align*}
a_1^2+\dots+a_9^2 &\leq 1/3+3\Delta^2+5(1/3-|\delta_1|/2)^2+(1/3-\delta_9)^2 \\
&= 1-5|\delta_1|/3-2\delta_9/3+5\delta_1^2/4+\delta_9^2+3\Delta^2 \\
&\leq 1-2\Delta/3+21\Delta^2/4.
\end{align*}

It follows that $\Var(Z_9)=1-a_1^2-\dots-a_9^2 \geq 2\Delta/3-21\Delta^2/4$.

Let $s=1$, $\gamma=9\Delta$, and $d=2/3$. Using the notations of \Cref{lem:LargeVarianceLeftNew} with $Y=Y_{9}$, we have
\begin{align*}
    p_{-1}&=\binom{9}{5}/2^9=63/256, \\
    p_0&=\binom{9}{6}/2^9=21/128, \\
    p_1&=\binom{9}{7}/2^9=9/128, \\
    t=\sum_{r=2}^\infty p_r=p_2+p_3&=\biggl[\binom{9}{8}+\binom{9}{9}\biggr]/2^9=5/256.
\end{align*}

Note that $p_{-1}\geq t$ and $p_0\geq p_1$.
Moreover, as $\Delta\leq 0.0009$ by \Cref{claim:case-c-claim-1}, $8.17\gamma^2=661.77\Delta^2$, and $2\Delta/3-21\Delta^2/4 \geq 661.77\Delta^2$ holds for all $\Delta\leq 0.000999$, we may apply \Cref{lem:LargeVarianceLeftNew} and conclude that
\begin{align*}
    \Prob(X \geq 1) &\geq p_0/4+3p_1/4+t
    = 29/256
    > 7/64,
\end{align*}
as desired.

\subsection{Case D}
\label{subsec:case-D}

Suppose for contradiction that there exists a counterexample to \Cref{thm:main-thm} with $a_i=1/4+\delta_i$ and $|\delta_i| \leq \delta \leq 0.03$ for $1 \leq i \leq 16$.

\begin{Claim}
\label{claim:case-d-claim-1}
    $a_1+a_2+a_3+a_4 < 1$, and $\delta\leq 0.00045$.
\end{Claim}
\begin{proof}
    It suffices to verify that the assumptions of \Cref{lem:small-sum} do hold for $c_i = 1/4$ for $1 \leq i \leq 16$, $I=\set{1, 2, 3, 4}$, $\lambda_i = 1$ for all $i\in I$, and $s=1$.

    Indeed, $\delta< 2(1/4)/16$, as required.
    We also have $|S|=\binom{16}{6}=8008$, $|R| = \binom{12}{6} = 924$, and $|T|=\binom{16}{5}+\dots+\binom{16}{0} = 6885$. 
    Thus $\min(|S|-|T|,|R|)/2 + |T| = 7374 \geq 2^{16}(7/64)$,
    Therefore by \Cref{lem:small-sum} we have $a_1+a_2+a_3+a_4 < 1$. 
    We may then feed this bound into the automated computation (see \Cref{appendix:simulation} \cref{sim:0D}) to find that we may strengthen our assumption on the size of $\delta$ to $\delta \leq 0.00045$.
\end{proof}

Let $\Delta=\max \{|\delta_1|, \dots, |\delta_{16}|\}$. Note that $\Delta=\max\{\delta_1, -\delta_{16} \}$. By \Cref{lem:IneqSumSquares}, we have
\begin{align*}
a_1^2+a_2^2+a_3^2+a_4^2 &\leq 1/4+4\Delta^2
\end{align*}

A computation of the same spirit of Case B and C will enable us to conclude in this case too, but we will need to be more careful in the computations of the constants so that we can apply \Cref{lem:LargeVarianceLeftNew} with the bound on $\delta$ from \Cref{claim:case-d-claim-1}, which we cannot improve without an unreasonably long computation.
To this end we split into two further subcases.

\subsubsection{Subcase 1: \texorpdfstring{$\Delta=\delta_1$}{Delta = delta\_1} or \texorpdfstring{$a_{14} \leq 1/4 -|\delta_{16}|/2$}{a\_14 <= 1/4 - |delta\_16|/2}}

First suppose that $\Delta=\delta_1$. Then by \Cref{claim:case-d-claim-1}, $a_4 \leq 1/4-|\delta_1|/3$. Hence
\begin{align*}
a_1^2+\dots+a_{16}^2 &\leq 1/4+4\Delta^2+12(1/4-|\delta_1|/3)^2 \\
&= 1-2|\delta_1|+4\delta_1^2/3+4\Delta^2 \\
&\leq 1-2\Delta+16\Delta^2/3.
\end{align*}

Hence $\Var(Z_{16})=1-a_1^2-\dots-a_{16}^2 \geq 2\Delta-16\Delta^2/3$ in this case.
If instead $\Delta=-\delta_{16}$, then by assumption we know that $a_{14} \leq 1/4 -|\delta_{16}|/2$.
Therefore
\begin{align*}
a_1^2+\dots+a_{16}^2 &\leq 1/4+4\Delta^2+9(1/4-|\delta_1|/3)^2+2(1/4 - |\delta_{16}|/2)^2+(1/4 - |\delta_{16}|)^2 \\
&= 1 -\delta_{16} + 4\Delta^2+3\delta_{16}^2/2 +\delta_1^2-3|\delta_1|/2 \\
&\leq  1-\Delta+11\Delta^2/2.
\end{align*}
Thus in this case $\Var(Z_{16})\geq \Delta-11\Delta^2/2$.
As $2\Delta-40\Delta^2/3\geq \Delta-11\Delta^2/2$ holds for all $\Delta \leq 6/47$, it follows that $\Var(Z_{16})\geq \Delta-11\Delta^2/2$.

We now seek to apply \Cref{lem:LargeVarianceLeftNew}. Let $s=1$, $\gamma=16\Delta$, and $d=1/2$. Using the notations of \Cref{lem:LargeVarianceLeftNew} with $Y=Y_{16}$, we have
\begin{align*}
    p_{-1}&=\binom{16}{9}/2^{16}=715/4096, \\
    p_0&=\binom{16}{10}/2^{16}=1001/8192, \\
    p_1&=\binom{16}{11}/2^{16}=273/4096, \\
    t=\sum_{r=2}^\infty p_r&=\biggl[\binom{16}{12}+\binom{16}{13}+\binom{16}{14}+\binom{16}{15}+\binom{16}{16}\biggr]/2^{16}=2517/65536.
\end{align*}

Note that $p_{-1}\geq t$ and $p_0\geq p_1$.
Moreover, as $\Delta\leq 0.00045$ by \Cref{claim:case-d-claim-1}, $8.17\gamma^2=2091.52\Delta^2$, and $\Delta-11\Delta^2/2 \geq 2091.52\Delta^2$ holds for all $\Delta\leq 0.000476$, we may apply \Cref{lem:LargeVarianceLeftNew} and conclude that
\begin{align*}
    \Prob(X \geq 1) &\geq p_0/4+3p_1/4+t
    = 7795/65536
    > 7/64,
\end{align*}
as required.

\subsubsection{Subcase 2: \texorpdfstring{$a_{14} > 1/4 -|\delta_{16}|/2$}{a\_14 > 1/4 - |delta\_16|/2}}
If $\Delta=\delta_1$ then we may proceed as in the previous subcase, so assume that $\Delta=-\delta_{16}$. 
In this case we have the following bound
\begin{align*}
a_1^2+\dots+a_{16}^2 &\leq 1/4+4\Delta^2+11(1/4)^2+(1/4 - |\delta_{16}|)^2 \\
&\leq  1-\Delta/2+5\Delta^2.
\end{align*}
It therefore follows that $\Var(Z_{16})=1-a_1^2-\dots-a_{16}^2 \geq \Delta/2-5\Delta^2$.

We now seek to apply \Cref{lem:LargeVarianceLeftNew}.  Let $s=1$, $\gamma=21\Delta/2$, and $d=1/2$, and note that we may indeed use $\gamma=21\Delta/2$ rather than $16\Delta$ as we have $\delta_4,\dotsc,\delta_{14}\in (-\Delta/2, 0]$.
We therefore see that $p_{-1},p_0,p_1,t$ are the same as in the previous subcase.
Furthermore, $8.17\gamma^2=988.57\Delta^2$, and $\Delta/2-5\Delta^2\geq 988.57\Delta^2$ holds for all $\Delta\leq 0.0005$, so \Cref{lem:LargeVarianceLeftNew} applies and we may conclude as in the previous subcase.

\subsection{Case E}
\label{subsec:case-E}

Suppose for contradiction that there exists a counterexample to \Cref{thm:main-thm} with $a_1=2/3+\delta_1$ and $a_i=1/3+\delta_i$ for $2 \leq i \leq 6$ such that $|\delta_i| \leq \delta \leq 0.002$ for each $1 \leq i \leq 6$. Let $\Delta=\max\{|\delta_1|, \dots, |\delta_6|\}$ and note that $\Delta= \max\{ |\delta_1|, -\delta_6 \}$ as $a_1+a_2 <1$ by \Cref{lem:a1+a2<1}, and $a_3+a_4+a_5 < 1$ by \Cref{lem:sum_a3a4a5}. Furthermore, as $a_1+a_2 <1$, we have
\begin{align*}
a_1^2+a_2^2 &< (2/3+ \delta_1)^2 + (1/3- \delta_1)^2 \\
&= 5/9 + 2\delta_1/3+2\delta_1^2 \\
&\leq 5/9 + 2\delta_1/3+2\Delta^2.
\end{align*}

As $a_3+a_4+a_5 < 1$, it follows by \Cref{lem:IneqSumSquares} that $a_3^2+a_4^2+a_5^2 \leq 1/3 +3\Delta^2$. Hence
\begin{align}
\label{eq:case-e-sum-sq-1}
a_1^2+a_2^2+a_3^2+a_4^2+a_5^2+a_6^2 &\leq 5/9 + 2\delta_1/3+2\Delta^2+1/3 +3\Delta^2+(1/3+\delta_6)^2 \nonumber \\
&= 1 +2\delta_1/3+2\delta_6/3 +6\Delta^2.
\end{align}

On the other hand, $a_2 \leq 1/3-\delta_1$, and consequently
\begin{align}
\label{eq:case-e-sum-sq-2}
a_1^2+a_2^2+a_3^2+a_4^2+a_5^2+a_6^2 &\leq (2/3+\delta_1)^2+4(1/3-\delta_1)^2+(1/3+\delta_6)^2 \nonumber \\
&= 1 -8\delta_1/3 +2\delta_6/3 +7\Delta^2.
\end{align}

By applying inequality \eqref{eq:case-e-sum-sq-1} if $\delta_1 < 0$ and \eqref{eq:case-e-sum-sq-2} if $\delta_1 \geq 0$, we find that $\Var(Z_6) \geq 2\Delta/3 - 7\Delta^2$.

Let $s=1$, $k=6$, $\gamma=6\Delta$, and $d=2/3$. Using the notations of \Cref{lem:LargeVarianceLeftNew} with $Y=Y_6$, we have
\begin{align*}
    p_{-1}&=\biggl[\binom{5}{2}+\binom{5}{4}\biggr]/2^6=15/64, \\
    p_0&=\biggl[\binom{5}{3}+\binom{5}{5}\biggr]/2^6=11/64, \\
    p_1&=\binom{5}{4}/2^6=5/64, \\
    t=\sum_{r=2}^\infty p_r=p_2&=\binom{5}{5}/2^6=1/64.
\end{align*}
Note that $p_{-1}\geq t$ and $p_0\geq p_1$.
Moreover, as $\Delta\leq 0.002$ (by \Cref{appendix:simulation} \cref{sim:0E}), $8.17\gamma^2=294.12\Delta^2$, and $2\Delta/3-7\Delta^2 \geq 294.12\Delta^2$ holds for all $\Delta\leq 0.0022$, we may apply \Cref{lem:LargeVarianceLeftNew} and conclude that
\begin{align*}
    \Prob(X \geq 1) &\geq p_0/4+3p_1/4+t
    = 15/128
    > 7/64,
\end{align*}
as desired.

\subsection{Case F}

Suppose for contradiction that there exists a counterexample to \Cref{thm:main-thm} with $a_1=1/2+\delta_1$ and $a_i=1/4+\delta_i$ for $2 \leq i \leq 13$ such that $|\delta_i| \leq \delta \leq 0.03$ for each $1 \leq i \leq 13$. \\
We will apply \Cref{lem:small-sum} twice.
In both applications, we have $c_1=1/2$ and $c_i=1/4$ for $2 \leq i \leq 13$, $s=1$, and $\lambda_i=1$ for all $i\in I$.
We also have $\delta\leq 0.03 < 2(1/4)/13$ as required by the lemma.
Furthermore, $|S|=\binom{12}{7}+\binom{12}{9}=1012$, and $|T|=\binom{12}{8}+\binom{12}{9}+2(\binom{12}{10}+\binom{12}{11}+\binom{12}{12}) = 873$.
The set $I$ and value of $|R|$, however, depend on the case.

\begin{Claim}
\label{claim:case-f-claim-1}
    We have $a_1+a_2+a_3 < 1$.
\end{Claim}

\begin{proof}
    In this case, $I=\set{1, 2, 3}$, and $|R| = \binom{10}{5} = 252$, and $\min(|S|-|T|,|R|)/2 + |T| = 942.5 \geq 2^{13}(7/64)$.
    Therefore by \Cref{lem:small-sum} we have $a_1+a_2+a_3 < 1$.
\end{proof}

\begin{Claim}
\label{claim:case-f-claim-2}
    We have $a_2+a_3+a_4+a_5 < 1$ and $\delta < 0.00035$.
\end{Claim}

\begin{proof}
    In this case, $I = \set{2, 3, 4, 5}$, and $|R| = 2\binom{8}{5} = 112$, and $\min(|S|-|T|,|R|)/2 + |T| = 929 \geq 2^{13}(7/64)$.
    Therefore by \Cref{lem:small-sum} we have $a_2+a_3+a_4+a_5 < 1$.
    Putting this bound and the bound from \Cref{claim:case-f-claim-1} into the automated computation \cite{simulation} (see \Cref{appendix:simulation} \cref{sim:0F}), we find that $\delta < 0.00035$, as required.
\end{proof}

Let $\Delta= \{ |\delta_1|, \dots, |\delta_{13}|\}$. Note that $\Delta= \max \{ |\delta_1|, |\delta_2|, |\delta_{13}| \}$.
If $\delta_1 \geq 0$, let $\Delta_{1,2}=\max \{|\delta_1|, |\delta_2| \}$ then $a_3 \leq 1/4-\Delta_{1,2}/2$ by \Cref{claim:case-f-claim-1}. Therefore
\begin{align}
\label{eq:case-f-sum-sq-1}
    a_1^2+\dots+a_{13}^2 &\leq (1/2+\Delta_{1,2})^2+(1/4+\Delta_{1,2})^2+10(1/4-\Delta_{1,2}/2)^2+(1/4+\delta_{13})^2 \nonumber \\
    &= 1 -\Delta_{1,2}+9\Delta_{1,2}^2/2+\delta_{13}/2+\delta_{13}^2 \nonumber \\
    &\leq 1-\Delta/2+11\Delta^2/2.
\end{align}

If $\delta_1 \leq 0$, using \Cref{lem:IneqSumSquares} with $a_2+a_3+a_4+a_5<1$ and the fact that $a_5 \leq 1/4-|\delta_2|/3$ by \Cref{claim:case-f-claim-2}, we get 
\begin{align}
\label{eq:case-f-sum-sq-2}
    a_1^2+\dots+a_{13}^2 &\leq (1/2)^2+1/4+4\Delta^2+7(1/4-|\delta_2|/3)^2+(1/4+\delta_{13})^2 \nonumber \\
    &= 1 - 7|\delta_2|/6+\delta_{13}/2+4\Delta^2+7\delta_2^2/9+\delta_{13}^2 \nonumber \\
    &\leq 1 -\Delta/2+52\Delta^2/9.
\end{align}

From inequalities \eqref{eq:case-f-sum-sq-1} and \eqref{eq:case-f-sum-sq-2}, we get $\Var(Z_{13})\geq \Delta/2 - 52\Delta^2/9$.
Let $s=1$, $\gamma=13\Delta$, and $d=1/2$. Using the notations from \Cref{lem:LargeVarianceLeftNew} with $Y=Y_{13}$, we have
\begin{align*}
    p_{-1}&=\biggl[\binom{12}{6}+\binom{12}{8}\biggr]/2^{13}=1419/8192, \\
    p_0&=\biggl[\binom{12}{7}+\binom{12}{9}\biggr]/2^{13}=253/2048, \\
    p_1&=\biggl[\binom{12}{8}+\binom{12}{10}\biggr]/2^{13}=561/8192, \\
    t=\sum_{r=2}^\infty p_r=p_2+p_3+p_4+p_5&=\biggl[\binom{12}{9}+\binom{12}{10}+2\binom{12}{11}+2\binom{12}{12}\biggr]/2^{13}=39/1024.
\end{align*}

Note that $p_{-1}\geq t$ and $p_0\geq p_1$.
Moreover, as $\Delta\leq 0.00035$ by \Cref{claim:case-f-claim-2}, $8.17\gamma^2=1380.73\Delta^2$, and $\Delta/2-52\Delta^2/9 \geq 1380.73\Delta^2$ holds for all $\Delta\leq 0.00036$, we may apply \Cref{lem:LargeVarianceLeftNew} and conclude that
\begin{align*}
    \Prob(X \geq 1) &\geq p_0/4+3p_1/4+t
    = 3943/32768
    > 7/64,
\end{align*}
as desired.

\subsection{Case G}

Suppose for contradiction that there exists a counterexample to \Cref{thm:main-thm} with $a_i=1/2+\delta_i$ for $1 \leq i \leq 2$ and $a_i=1/4+\delta_i$ for $3 \leq i \leq 10$ such that $|\delta_i| \leq \delta \leq 0.05$ for each $1 \leq i \leq 10$.

We will apply \Cref{lem:small-sum} three times. In all three applications, we set $c_1=c_2=1/2$, and $c_i=1/4$ for $3\leq i \leq 10$, and $s=1$. We have that $|S|=\binom{8}{4}+2\binom{8}{6}+\binom{8}{8}=127$, and $|T|=\binom{8}{5}+\binom{8}{6}+3\binom{8}{7}+3\binom{8}{8} = 111$.
We also have that $\delta\leq 0.005 < 2(1/4)/10$, and so it remains in each claim to compute the value of $|R|$, and check that the required inequality holds.

\begin{Claim}
\label{claim:case-g-claim-1}
    We have $a_2+a_3+a_4 < 1$.
\end{Claim}
\begin{proof}
    We set $I=\set{2, 3, 4}$ and $\lambda_i=1$ for each $i\in I$. We have that $|R| = 2\binom{6}{4} = 30$, and so $\min(|S|-|T|,|R|)/2 + |T| = 119 \geq 2^{10}(7/64)$. The result then follows from \Cref{lem:small-sum}.
\end{proof}

\begin{Claim}
\label{claim:case-g-claim-2}
    We have $a_3+a_4+a_5+a_6 < 1$ and $\delta < 0.00045$.
\end{Claim}
\begin{proof}
    We set $I=\set{3, 4, 5, 6}$ and $\lambda_i = 1$ for each $i\in I$. We have that $|R| = \binom{4}{0} + 2\binom{4}{2} + \binom{4}{4} = 14$, and so $\min(|S|-|T|,|R|)/2 + |T| = 118 \geq 2^{10}(7/64)$. The result follows from \Cref{lem:small-sum}.
    Putting this bound and the bound from \Cref{claim:case-g-claim-1} into the automated computation (see \Cref{appendix:simulation} \cref{sim:0G}), we find that $\delta < 0.00045$, as required.
\end{proof}

\begin{Claim}
    We have $a_1-a_2+a_3+a_4+a_5+a_6 < 1$.
\end{Claim}

\begin{proof}
    We set $I=\set{1,2,3,4,5,6}$, $\lambda_2=-1$, and $\lambda_i=1$ for each $i\in I\setminus \set{2}$. We have that $|R| = \binom{4}{2} = 6$, and so $\min(|S|-|T|,|R|)/2 + |T| = 114 \geq 2^{10}(7/64)$.
    Therefore by \Cref{lem:small-sum} we have $a_1-a_2+a_3+a_4+a_5+a_6 < 1$.
\end{proof}

By \Cref{lem:IneqSumSquares}, $a_1^2+a_2^2 \leq 1/2 + 2\Delta^2$. By \Cref{lem:IneqSumSquares} with $b=(1-\delta_1+\delta_2)/4$ and $\epsilon=3\Delta/2$, we also have 
\begin{align*}
    a_3^2+a_4^2+a_5^2+a_6^2 &\leq (1-\delta_1+\delta_2)^2/4+6\Delta^2 \\
    &\leq 1/4-\delta_1/2+\delta_2/2-\delta_1\delta_2/2+13\Delta^2/2.
\end{align*}

Moreover, $a_3+3a_6 \leq a_3+a_4+a_5+a_6 \leq a_1-a_2+a_3+a_4+a_5+a_6 < 1$, so $a_6 \leq 1/4- |\delta_3|/3$. It follows that
\begin{align*}
    a_1^2+\dots+a_{10}^2 &\leq 1/2 + 2\Delta^2+1/4-\delta_1/2+\delta_2/2-\delta_1\delta_2/2+13\Delta^2/2+4(1/4-|\delta_3|/3)^2 \\
    &\leq 1- \Delta/2 +85\Delta^2/9.
\end{align*}

Hence $\Var(Z_{10}) \geq \Delta/2 - 85\Delta^2/9$.
Let $s=1 $, $\gamma=10\Delta$, and $d=1/2$. Using the notations of \Cref{lem:LargeVarianceLeftNew} with $Y=Y_{10}$, we have
\begin{align*}
    p_{-1}&=\biggl[\binom{8}{3}+2\binom{8}{5}+\binom{8}{7}\biggr]/2^{10}=11/64, \\
    p_0&=\biggl[\binom{8}{4}+2\binom{8}{6}+\binom{8}{8}\biggr]/2^{10}=127/1024, \\
    p_1&=\biggl[\binom{8}{5}+2\binom{8}{7}\biggr]/2^{10}=9/128, \\
    t=\sum_{r=2}^\infty p_r=p_2+p_3+p_4&=\biggl[\binom{8}{6}+\binom{8}{7}+3\binom{8}{8}\biggr]/2^{10}=39/1024.
\end{align*}

Note that $p_{-1}\geq t$ and $p_0\geq p_1$.
Moreover, as $\Delta\leq 0.00045$ from \Cref{claim:case-g-claim-2}, $8.17\gamma^2=817\Delta^2$, and $\Delta/2-85\Delta^2/9 \geq 817\Delta^2$ holds for all $\Delta\leq 0.000605$, we may apply \Cref{lem:LargeVarianceLeftNew} and conclude that
\begin{align*}
    \Prob(X \geq 1) &\geq p_0/4+3p_1/4+t
    = 499/4096
    > 7/64,
\end{align*}
as desired.

\subsection{Case H}
Note that this case is not technically a hard case as we defined it, since for $a_1=a_2=a_3=1/2$ and $a_4=a_5=a_6=a_7=1/4$, we have $\Prob(X > 1)=7/64$. However, we can deal with this case in the same manner as the previous cases.

Suppose for contradiction that there exists a counterexample to \Cref{thm:main-thm} with $a_i=1/2+\delta_i$ for $1 \leq i \leq 3$ and $a_i=1/4+\delta_i$ for $4 \leq i \leq 7$ such that $|\delta_i| \leq \delta \leq 0.005$ for each $1 \leq i \leq 7$.
Letting $\Delta=\{|\delta_1|, \dots, |\delta_7|\}$ and $k=7$, we have that $\Delta < 0.005$ by \Cref{appendix:simulation} \cref{sim:0B}. We arrive at the following by simple counting, nothing that the intervals are pairwise disjoint.
\begin{align*}
    p_0&=\Prob(Y_7\in [1-7\Delta,1+7\Delta])=\biggl[\binom{3}{0}\binom{4}{3}+\binom{3}{1}\binom{4}{1}\biggr]/2^{7}=16/128 \\
    p_1&=\Prob(Y_7\in [3/2-7\Delta,3/2+7\Delta])=\biggl[\binom{3}{0}\binom{4}{2}+\binom{3}{1}\binom{4}{4}\biggr]/2^{7}=9/128 \\
    p_2&=\Prob(Y_7\in [2-7\Delta,2+7\Delta])=\binom{3}{0}\binom{4}{3}/2^{7}=4/128 \\
    p_3&=\Prob(Y_7\in [5/2-7\Delta,5/2+7\Delta])=\binom{3}{0}\binom{4}{4}/2^{7}=1/128.
\end{align*}

But then we can see that, letting $q=\Prob(Z_7\geq 7\Delta)$,
\begin{align*}
    \Prob(X\geq 1)&\geq qp_0+(1-q)(p_1+p_2+p_3) \\
    &= (16q+(1-q)14)/128\\
    &\geq 14/128 \\
    &= 7/64,
\end{align*}
as desired.

\section{Lowther's conjecture}
\label{sec:lowther-conjecture}

In this section we prove the following results.

\begin{theorem}
\label{thm:sqrt-7-thm}
For any Rademacher sum $X = \sum_{i=1}^n a_i\eps_i\in\RadSums$, we have

$$\mathbb{P}(|X| \geq 1/\sqrt{7}) \geq 1/2.$$
\end{theorem}

\begin{theorem}
\label{thm:sqrt-5-thm}
    For any Rademacher sum $X=\sum_{i=1}^n a_i\eps_i\in\RadSums$, we have
    $$\Prob(|X|\geq 1/\sqrt{5})\geq 29/64.$$
\end{theorem}

\begin{theorem}
\label{thm:sqrt-3-thm}
    For any Rademacher sum $X=\sum_{i=1}^n a_i\eps_i\in\RadSums$, we have
    $$\Prob(|X|\geq 1/\sqrt{3})\geq 3/8.$$
\end{theorem}

\begin{theorem}
\label{thm:sqrt-6-thm}
For any Rademacher sum $X = \sum_{i=1}^n a_i\eps_i \in\RadSums$, we have

$$\mathbb{P}(|X| \geq 2/\sqrt{6}) \geq 1/4.$$
\end{theorem}

Note that the above Theorems imply \Cref{thm:ExplicitInfFunction}.

\begin{proof}[Proof of \Cref{thm:ExplicitInfFunction} assuming Theorems \ref{thm:sqrt-7-thm}, \ref{thm:sqrt-5-thm}, \ref{thm:sqrt-3-thm}, \ref{thm:sqrt-6-thm} and \ref{thm:main-thm}]
Theorems \ref{thm:sqrt-7-thm}, \ref{thm:sqrt-5-thm}, \ref{thm:sqrt-3-thm}, \ref{thm:sqrt-6-thm} and \ref{thm:main-thm} already give us a lower bound, that is,

\begin{itemize}
    \item $f(y) \geq 1/2$ for $y \in (0,1/\sqrt{7}]$,
    \item $f(y) \geq 29/64$ for $y \in (1/\sqrt{7},1/\sqrt{5}]$,
    \item $f(y) \geq 3/8$ for $y \in (1/\sqrt{5},1/\sqrt{3}]$,
    \item $f(y) \geq 1/4$ for $y \in (1/\sqrt{3},2/\sqrt{6}]$,
    \item $f(y) \geq 7/32$ for $y \in (2/\sqrt{6},1]$.
\end{itemize}

To see that all those bounds are in fact equalities, we remark that

\begin{itemize}
    \item $(a_1, a_2)=(1/\sqrt{2},1/\sqrt{2})$ shows that $f(y) \leq 1/2$ for $y \in (0,1/\sqrt{7}]$,
    \item $(a_1, \dots, a_7)=(1/\sqrt{7}, \dots, 1/\sqrt{7})$ shows that $f(y) \leq 29/64$ for $y \in (1/\sqrt{7},1/\sqrt{5}]$,
    \item $(a_1, \dots, a_5)=(1/\sqrt{5}, \dots, 1/\sqrt{5})$ shows that $f(y) \leq 3/8$ for $y \in (1/\sqrt{5},1/\sqrt{3}]$,
    \item $(a_1, \dots, a_3)=(1/\sqrt{3}, \dots, 1/\sqrt{3})$ shows that $f(y) \leq 1/4$ for $y \in (1/\sqrt{3},2/\sqrt{6}]$,
    \item $(a_1, \dots, a_6)=(1/\sqrt{6}, \dots, 1/\sqrt{6})$ shows that $f(y) \leq 7/32$ for $y \in (2/\sqrt{6},1]$,
    \item $a_1=1$ shows that $f(y) \leq 0$ for $y \in (1,\infty)$
\end{itemize}
    
Moreover, it is clear that $f(0)=1$. This finishes the proof of \Cref{thm:ExplicitInfFunction}.
\end{proof}

We now proceed to prove Theorems \ref{thm:sqrt-7-thm}, \ref{thm:sqrt-5-thm}, \ref{thm:sqrt-3-thm} and \ref{thm:sqrt-6-thm}, by methods mostly similar to those used in \Cref{sec:lowther-conjecture}.

In the cases $y = 1/\sqrt{7}$ and $y=2/\sqrt{6}$, (\Cref{subsec:sqrt-7-proof} and \Cref{subsec:2-sqrt-6-proof} respectively) we can prove an upper bound on the $a_i$'s which rules out all equality cases, leaving us with only a single neighbourhood of a hard case to consider, which can be resolved using arguments similar to those used in the previous section.
In the other two cases, namely $y=1/\sqrt{5}$ and $y=1/\sqrt{3}$, (\Cref{subsec:sqrt-5-proof} and \Cref{subsec:sqrt-3-proof} respectively) we cannot rule out the equality cases via a simple upper bound on the $a_i$'s. 
For this reason we split both of these cases into two subcases, one of which will be a neighbourhood of a hard case, and the other will contain the equality cases.

\subsection{Proof of Theorem \ref{thm:sqrt-7-thm}}
\label{subsec:sqrt-7-proof}

This case is the easiest of all of the cases we deal with, and can be resolved with a short preliminary result and a single application of \Cref{lem:LargeVarianceLeftNew}. Assume for contradiction that $X =\sum_{i=1}^n a_i\epsilon_i \in \mathcal{X}$ satisfies $\mathbb{P}(X \geq 1/\sqrt{7}) < 1/4$.

\begin{Claim}
\label{claim:sqrt-7-claim-1}
    $a_1 < 1/\sqrt{7}$, and $a_i=1/\sqrt{7}+\delta_i$ for $1\leq i \leq 7$, with each $\delta_i\in ( -0.001,0)$.
\end{Claim}
\begin{proof}
    If $a_1 \geq 1/\sqrt{7}$, then applying \Cref{obs:observation1DK} gives us that $\mathbb{P}(|X| \geq 1/\sqrt{7}) \geq 1/2$.
    Therefore we may assume that $a_1 < 1/\sqrt{7}$.
    
    Using the bound $a_1 < 1/\sqrt{7}$ in the automated computation (see \Cref{appendix:simulation} \cref{sim:1}), we find that we may assume further that $a_i=1/\sqrt{7}+\delta_i$ for $1\leq i \leq 7$ such that each $|\delta_i| < 0.001$, as required. 
\end{proof}

We now derive a contradiction via \Cref{lem:LargeVarianceLeftNew}. Let $\Delta=\max\set{|\delta_1|,\dotsc,|\delta_7|}=-\delta_7$.
Let $s=1/\sqrt{7}$, $\gamma=7\Delta$, and $d=2/\sqrt{7}$.
Using the notations of \Cref{lem:LargeVarianceLeftNew} with $Y=Y_7$, we have
\begin{align*}
    p_{-1}&=\binom{7}{3}/2^7 = 35/128, \\
    p_0&=\binom{7}{4}/2^7 = 35/128,\\
    p_1&=\binom{7}{5}/2^7 = 21/128,\\
    t=\sum_{t=2}^\infty=p_2+p_3&=\biggl[\binom{7}{6} + \binom{7}{7}\biggr]/2^7 = 8/128.
\end{align*}
Next, $\Var(Z_7)\geq 1 - 6(1/\sqrt{7})^2-(1/\sqrt{7}-\Delta)^2 = 2\Delta/\sqrt{7}-\Delta^2$ and $8.17\gamma^2 = 400.33\Delta^2$, and so $\Var(Z_7)\geq 8.17\gamma^2$ holds whenever $\Delta\leq 0.00188$.
As \Cref{claim:sqrt-7-claim-1} gives $\Delta\leq 0.001$, and $p_{-1}\geq t$ and $p_0\geq p_1$, we may apply \Cref{lem:LargeVarianceLeftNew} with $s,\gamma,d$ as above to derive
\begin{align*}
    \Prob(X\geq 1/\sqrt{7})\geq p_0/4+3p_1/4+t = 65/256 > 1/4,
\end{align*}
as required.

\subsection{Proof of Theorem \ref{thm:sqrt-5-thm}}
\label{subsec:sqrt-5-proof}

Assume for contradiction that $X =\sum_{i=1}^n a_i\epsilon_i \in \mathcal{X}$ satisfies $\Prob(X\geq 1/\sqrt{5}) < 29/128$. As discussed at the start of \Cref{sec:lowther-conjecture}, we split into two subcases, one containing the equality cases, and the other a neighbourhood of a hard case.
We make this split depending on the value of $a_1+a_2+a_3+a_4+a_5$.

\subsubsection{Case: \texorpdfstring{$a_1+a_2+a_3+a_4+a_5\leq 2.1$}{a1+a2+a3+a4+a5 <= 2.1}}
\label{subsubsec:sqrt-5-small-sum}

This deals with the equality cases.
We wish to apply \Cref{lem:pairing}, and to that end we set
\begin{align*}
    A=\set{\zeta\in\set{-1,+1}^7\st\sum_{i=1}^7\zeta_i \geq 3}, \\
    B=\set{\zeta\in\set{-1,+1}^7\st\sum_{i=1}^7\zeta_i = 1},
\end{align*}
and note that $|B| = \binom{7}{3} = 35 > 29 = \binom{7}{2}+\binom{7}{1}+\binom{7}{0}=|A|$.

\begin{Claim}
    We have the following two inequalities.
    \begin{subequations}
    \label{eq:sqrt-5-sum-lower-bound}
        \begin{align}
            -a_1-a_2+a_3+a_4+a_5+a_6+a_7 &\geq 0.95. \label{eq:sqrt-5-sum-lower-bound-1}\\
            -a_1-a_2-a_3+a_4+a_5+a_6+a_7 &\geq 0.175. \label{eq:sqrt-5-sum-lower-bound-2}
        \end{align}
    \end{subequations}
\end{Claim}
\begin{proof}
    By automated computation; see \Cref{appendix:simulation} \cref{sim:2}.
\end{proof}

Now, we may note that $0.175 + 0.95 = 1.125 > 2/\sqrt{5}$, and the bounds in equations \eqref{eq:sqrt-5-sum-lower-bound-1} and \eqref{eq:sqrt-5-sum-lower-bound-2} are also lower bounds for any sum with signs taken from $A$ and $B$ respectively. 
Therefore \Cref{lem:pairing} applies with $s=1/\sqrt{5}$, $k=7$, and an arbitrary injection $\phi\st A\to B$.
Hence, eliminating $a_1,\dotsc,a_7$ by \Cref{obs:elimination}, we have that 
\begin{align*}
    \Prob(X\geq s)\geq 2^{-7}\sum_{\eps\in A\cup B}\Prob\Bigl[\sum_{i=1}^7\eps_i a_i + Z_7\geq s\Bigr]\geq 2^{-7}|A| = 29/128.
\end{align*}
This is exactly as required to conclude this case.

\subsubsection{Case: \texorpdfstring{$a_1+a_2+a_3+a_4+a_5 > 2.1$}{a1+a2+a3+a4+a5 > 2.1}}
\label{subsubsec:sqrt-5-big-sum}

Let $a_i=1/\sqrt{5}+ \delta_i$ for $1 \leq i \leq 5$.

\begin{Claim}
\label{claim:sqrt-5-subcase-2-claim-1}
    We have $a_1 < 1/\sqrt{5}$, and each $\delta_i\in ( -0.0025,0)$.
\end{Claim}

\begin{proof}
    If $a_1 \geq 1/\sqrt{5}$, then we get a contradiction by \Cref{obs:observation1DK}.\\
    Adding the bound $a_1 < 1/\sqrt{5}$ into the automated computation (see \Cref{appendix:simulation} \cref{sim:2}), we find that we may furthermore assume that each $\delta_i\in (-0.0025,0)$.
\end{proof}

Hence $\Delta=\max \{ |\delta_1|, \dots, |\delta_5| \} = - \delta_5$ and we have
\begin{align*}
    a_1^2+\dots+a_5^2 &\leq 4/5+(1/\sqrt{5}+\delta_5)^2 \\
    &= 1-2\Delta/\sqrt{5}+\Delta^2.
\end{align*}

We have $\Var(Z_5)\geq 2\Delta/\sqrt{5} - \Delta^2$.
Let $s=1/\sqrt{5} $, $\gamma=5\Delta$, and $d=2/\sqrt{5}$. Using the notations of \Cref{lem:LargeVarianceLeftNew} with $Y=Y_5$, we have
\begin{align*}
    p_{-1}&=\binom{5}{2}/2^{5}=5/16, \\
    p_0&=\binom{5}{3}/2^{5}=5/16, \\
    p_1&=\binom{5}{4}/2^{5}=5/32, \\
    t=\sum_{r=2}^\infty p_r=p_2&=\binom{5}{5}/2^{5}=1/32.
\end{align*}
Note that $p_{-1}\geq t$ and $p_0\geq p_1$.
Moreover, as $\Delta\leq 0.0025$ by \Cref{claim:sqrt-5-subcase-2-claim-1}, $8.17\gamma^2=204.25\Delta^2$, and $2\Delta/\sqrt{5}-\Delta^2 \geq 204.25\Delta^2$ holds for all $\Delta\leq 0.0043$, we may apply \Cref{lem:LargeVarianceLeftNew} and conclude that
\begin{align*}
    \Prob(X \geq 1) &\geq p_0/4+3p_1/4+t
    = 29/128,
\end{align*}
as desired.

\subsection{Proof of Theorem \ref{thm:sqrt-3-thm}}
\label{subsec:sqrt-3-proof}
Assume for contradiction that $X =\sum_{i=1}^n a_i\epsilon_i \in \mathcal{X}$ satisfies $\Prob(X\geq 1/\sqrt{3}) < 3/16$.
As in \Cref{subsec:sqrt-5-proof}, we split into two cases, but this time depending on the value of $a_1+a_2+a_3$.

\subsubsection{Case: \texorpdfstring{$a_1+a_2+a_3\leq 1.6$}{a1+a2+a3 <= 1.6}}
\label{subsubsec:sqrt-3-small-sum}

This deals with the equality cases, and proceeds analogously to \Cref{subsubsec:sqrt-5-small-sum}.
We define
\begin{align*}
    A=\set{\zeta\in\set{-1,+1}^5\st\sum_{i=1}^5\zeta_i \geq 3}, \\
    B=\set{\zeta\in\set{-1,+1}^5\st\sum_{i=1}^5\zeta_i = 1}, 
\end{align*}
and have the following claim.

\begin{Claim}
    We have the following two inequalities.
    \begin{subequations}
    \label{eq:sqrt-3-sum-lower-bounds}
        \begin{align}
            -a_1+a_2+a_3+a_4+a_5 \geq 1.04. \label{eq:sqrt-3-sum-lower-bound-1}\\
            -a_1-a_2+a_3+a_4+a_5 \geq 0.23. \label{eq:sqrt-3-sum-lower-bound-2}
        \end{align}
    \end{subequations}
\end{Claim}
\begin{proof}
    By automated computation. See \Cref{appendix:simulation} \cref{sim:3}.
\end{proof}

Note that inequalities \eqref{eq:sqrt-3-sum-lower-bound-1} and \eqref{eq:sqrt-3-sum-lower-bound-2} also give bounds on the sums of $a_1, \dots,a_5$ with signs taken from $A$ and $B$ respectively.
Note further that $0.23+ 1.04=1.27 > 2/\sqrt{3}$.
Thus if we take $\phi\st A\to B$ an arbitrary injection (noting that $|A|=6 < 10=|B|$), then \Cref{lem:pairing} applies with $s=1/\sqrt{3}$ and $k=5$, and we may conclude that $\Prob(X\geq 1/\sqrt{3})\geq 6/2^5 = 3/16$, as required.

\subsubsection{Case: \texorpdfstring{$a_1+a_2+a_3 > 1.6$}{a1+a2+a3 > 1.6}}
\label{subsubsec:sqrt-3-big-sum}

Let $a_i=1/\sqrt{3}+\delta_i$ for each $1 \leq i \leq 3$.

\begin{Claim}
\label{claim:sqrt-3-subcase-2-claim-1}
    $a_1 < 1/\sqrt{3}$, and each $\delta_i \in (-0.001,0)$.
\end{Claim}

\begin{proof}
    If $a_1 \geq 1/\sqrt{3}$, then we get a contradiction by \Cref{obs:observation1DK}. \\
    Adding the bound $a_1 < 1/\sqrt{3}$ into the automated computation (see \Cref{appendix:simulation} \cref{sim:3}), we find that we may furthermore assume that each $\delta_i\in (-0.001,0)$.
\end{proof}

Hence $\Delta=\max \{|\delta_1|, |\delta_2|, |\delta_3| \}=-\delta_3$. Furthermore
\begin{align*}
    a_1^2+a_2^2+a_3^2 &\leq (1/\sqrt{3})^2+(1/\sqrt{3})^2+(1/\sqrt{3}+\delta_3)^2 \\
    &= 1-2/\sqrt{3}\Delta+\Delta^2.
\end{align*}

so $\Var(Z_3) = 2/\Delta/\sqrt{3}-\Delta^2$.
Let $s=1/\sqrt{3}$, $\gamma=3\Delta$, and $d=2/\sqrt{3}$. Using the notations of \Cref{lem:LargeVarianceLeftNew} with $Y=Y_3$, we have
\begin{align*}
    p_{-1}&=\binom{3}{1}/2^3=3/8, \\
    p_0&=\binom{3}{2}/2^3=3/8, \\
    p_1&=\binom{3}{3}/2^3=1/8, \\
    t&=\sum_{r=2}^\infty p_r=0.
\end{align*}
Note that $p_{-1}\geq t$ and $p_0\geq p_1$.
Moreover, as $\Delta\leq 0.001$ by \Cref{claim:sqrt-3-subcase-2-claim-1}, $8.17\gamma^2=73.53\Delta^2$, and $2\Delta/\sqrt{3}-\Delta^2 \geq 73.53\Delta^2$ holds for all $\Delta\leq 0.015$, we may apply \Cref{lem:LargeVarianceLeftNew} and conclude that
\begin{align*}
    \Prob(X \geq 1) &\geq p_0/4+3p_1/4+t
    \geq 3/16,
\end{align*}
as desired.

\subsection{Proof of Theorem \ref{thm:sqrt-6-thm}}
\label{subsec:2-sqrt-6-proof}

Assume for contradiction that $X =\sum_{i=1}^n a_i\epsilon_i \in \mathcal{X}$ satisfies $\Prob(X\geq 2/\sqrt{6}) < 1/8$.

\begin{Claim}
\label{a_1+a_2 <2/sqrt{6}}
    We have $a_1+a_2 <2/\sqrt{6}$, and $a_i=1/\sqrt{6}+\delta_i$ with $|\delta_i|\leq 0.0025$ for each $1\leq i \leq 6$.
\end{Claim}

\begin{proof}
    First note that $a_1+a_2<2/\sqrt{6}$ follows immediately from \Cref{obs:observation1DK}.
    Adding this constraint to the automated computation (See \Cref{appendix:simulation} \cref{sim:4}), we find that $|\delta_i|\leq 0.0025$, as required.
\end{proof}

Let $\Delta=\max \{ |\delta_1|, \dots, |\delta_6| \}$ and note that $\Delta= \max \{ \delta_1, -\delta_6 \}$. 
Note further that $a_3 \leq 1/\sqrt{6}-|\delta_1|$ by \Cref{a_1+a_2 <2/sqrt{6}}. 
Using \Cref{lem:IneqSumSquares}, we have
\begin{align*}
    a_1^2+\dots+a_6^2 &\leq 1/6 + 2\Delta^2+3(1/\sqrt{6}-|\delta_1|)^2+(1/\sqrt{6}+\delta_6)^2 \\
    &= 1-\sqrt{6}|\delta_1|+2\delta_6/\sqrt{6}+ 3\delta_1^2+ \delta_6^2+ 2\Delta^2 \\
    &\leq 1-2\Delta/\sqrt{6}+6\Delta^2,
\end{align*}
so $\Var(Z_6)\geq 2\Delta/\sqrt{6} - 6\Delta^2$.
Let $s=2/\sqrt{6}$, $\gamma=6\Delta$, and $d=2/\sqrt{6}$. Using the notations of \Cref{lem:LargeVarianceLeftNew} with $Y=Y_6$, we have
\begin{align*}
    p_{-1}&=\binom{6}{3}/2^{6}=5/16, \\
    p_0&=\binom{6}{4}/2^{6}=15/64, \\
    p_1&=\binom{6}{5}/2^{6}=3/32, \\
    t=\sum_{r=2}^\infty p_r=p_2&=\binom{6}{6}/2^{6}=1/64.
\end{align*}
Moreover, as $\Delta\leq 0.0025$ by \Cref{a_1+a_2 <2/sqrt{6}}, $8.17\gamma^2=294.12\Delta^2$, and $2\Delta/\sqrt{6}-6\Delta^2 \geq 294.12\Delta^2$ holds for all $\Delta\leq 0.00272$, we may apply \Cref{lem:LargeVarianceLeftNew} and conclude that
\begin{align*}
    \Prob(X \geq 1) &\geq p_0/4+3p_1/4+t
    = 37/256
    > 1/8,
\end{align*}
as desired.

\section{Concluding remarks and further work}

Recall that $\RadSums$ is the class of all Rademacher sums with variance 1. Keller and Klein \cite{keller2022proof} asked the following.

\begin{Problem}
    Characterise the function $F$, defined for all $x \in \mathbb{R}$, by $F(x)= \sup_{X \in \RadSums} \mathbb{P}(X > x)$.
\end{Problem}

Note that $F(x)= 1- \inf_{X \in \RadSums} \mathbb{P}(X \geq -x)$, so \Cref{thm:ExplicitInfFunction} answers the question of Keller and Klein for all negative $x$. Pinelis \cite{pinelis2012asymptotically} gave asymptotic values of $F$, and some precise values of $F$ have been determined, but most of the behaviour of $F$ remains to be understood. It had been previously conjectured that for every $x$, the $\sup$ in the definition of $F(x)$ is attained for Rademacher sums having all of their weights equal. However, this conjecture was disproved by Pinelis \cite{pinelis2015supremum}.

While it seems that determining $F(x)$ for all large $x$ is beyond the reach of current techniques, it appears that the function $F$ stays moderately well-behaved for some values of $x$ greater than 1.
Indeed, Keller and Klein \cite[Theorem 1.2]{keller2022proof} showed that $F(x)=1/4$ for $x\in[1,\sqrt{2})$.
We make the following conjecture for larger values of $x$.

\begin{conjecture}
    For $x\in[\sqrt{2},2]$, $F(x)$ is characterised as follows.
    \begin{itemize}
        \item $F(x)=1/8$ for $x\in[\sqrt{2},\sqrt{3})$,
        \item $F(x)=1/16$ for $x\in[\sqrt{3},2)$,
        \item $F(2)=9/256$.
    \end{itemize}
\end{conjecture}

\bibliographystyle{abbrvnat}  
\renewcommand{\bibname}{Bibliography}
\bibliography{main}

\appendix

\section{Bounds proved by automated computation}
\label{appendix:simulation}
The bounds proved by the automated computation carried out by our computer program are stored in a number of files accessible with the code \cite{simulation}, in the directory \texttt{cases}.
Here we provide a list of all the bounds proved by automated computation which we use, and the files which contain them.
The filenames are designed to align with the structure of the paper.
For example, Case C of \Cref{sec:tomaszewski-counterpart} is contained within files beginning \texttt{0C}.

\begin{enumerate}
    \item \label{sim:1} The $s=1/\sqrt{7}$ case is dealt with in file \texttt{1}. The restriction that $a_1\leq 0.378$ is enforced, as it is proved in \Cref{claim:sqrt-7-claim-1} that $a_1\leq 1/\sqrt{7} < 0.378$ before any results from the automated computation are used. The automated computation then proves that $|a_i-1/\sqrt{7}|\leq 0.001$ for $1\leq i\leq 7$.
    \item \label{sim:2} The $s=1/\sqrt{5}$ case is dealt with in file \texttt{2}. We split into two subcases depending on whether the value of $a_1+\dots+a_5$ is greater than or less than $2.1$, and these subcases are further investigated in files \texttt{2A} and \texttt{2B} respectively. 
    In file \texttt{2A} it is shown that $|a_i-1/\sqrt{5}|\leq 0.0025$ for each $1\leq i\leq 5$.
    In file \texttt{2B} it is shown that $-a_1-a_2+a_3+\dots+a_7\geq 0.95$ and $-a_1-a_2-a_3+a_4+\dots+a_7\geq 0.175$.
    \item \label{sim:3} The $s=1/\sqrt{3}$ case is dealt with in file \texttt{3}. We split into two subcases depending on whether the value of $a_1+a_2+a_3$ is greater than or less than $1.6$, and these subcases are further investigated in files \texttt{3A} and \texttt{3B} respectively.
    In file \texttt{3A} it is shown that $|a_i-1/\sqrt{3}|\leq 0.001$ for each $1\leq i\leq 3$.
    In file \texttt{3B} it is shown that $-a_1+a_2+a_3+a_4+a_5\geq 1.04$ and $-a_1-a_2+a_3+a_4+a_5\geq 0.23$.
    \item \label{sim:4} The $s=2/\sqrt{6}$ case is dealt with in file \texttt{4}. The restriction that $a_1+a+2\leq 0.817$ is enforced, as it is proved in \Cref{a_1+a_2 <2/sqrt{6}} that $a_1+a+2\leq 2/\sqrt{6} < 0.817$ before any results from the automated computation are used. The automated computation then proves that $|a_i-1/\sqrt{6}|\leq 0.0025$ for $1\leq i\leq 6$.
    \item \label{sim:0} The $s=1$ case is split into subcases $A-L$ via instructions in file \texttt{0}. Subcases $I-L$ exist primarily for ease of computation, and could in theory have been combined into a single subcase from which a contradiction could be derived. The subcases and how they are resolved are as follows.
    \begin{enumerate}[label=5.\Alph*.]
        \item \label{sim:0A} $0.9\leq a_1$. \\
        File \texttt{0A} proves that $a_1\in (1 - \delta, 1]$ for $\delta=0.04$.
        \item \label{sim:0B} $a_1+a_2+a_3\geq 1.4$ and $a_1\leq 0.6$. \\
        Note that this case includes both cases B and H in the division into cases of \Cref{sec:lowther-conjecture}. 
        File \texttt{0B} deals with both cases; splitting on whether $a_4$ is more or less than $0.35$. 
        In both cases it demonstrates that each for each $a_i$ there is some $c_i\in\set{0, 0.25, 0.5}$ such that $|a_i-c_i|\leq 0.005$ for each $1\leq i \leq 8$.
        Noting that we must have $c_3=0.5$ we arrive at the two expected cases.
        \item \label{sim:0C} $0.9\leq a_1+a_2+a_3\leq 1.1$ and $a_1\leq 0.45$.\\ 
        This case is difficult to pin down, and so to speed up computation, file \texttt{0C} first splits into four subcases based on the value of $a_7$, and we then deal with each individually.
        \begin{enumerate}[label=5.C.\Alph*.]
            \item $a_7\geq 0.325$. \\
            File \texttt{0CA} proves that $|a_i-1/3|\leq 0.07$ for each $1\leq i \leq 9$, as this allows \Cref{claim:case-c-claim-1} to apply, and then file \texttt{0CAA} can improve upon these bounds by enforcing the extra assumption that $a_1+a_2+a_3 < 1$. 
            The resulting bounds are then fed into file \texttt{0CAAA} (for efficiency purposes) which proceeds to prove the required bound that $|a_i-1/3|\leq 0.0009$ for $1\leq i \leq 9$.
            \item $0.3\leq a_7\leq 0.325$. \\
            File \texttt{0CB} proves that $|a_i-1/3|\leq 0.07$ for each $1\leq i \leq 9$, as this allows \Cref{claim:case-c-claim-1} to apply, and then file \texttt{0CBA} can the enforce the extra assumption that $a_1+a_2+a_3 < 1$, and with this extra information, can derive a contradiction.
            \item $0.27\leq a_7\leq 0.3$. \\
            File \texttt{0CC} derives a contradiction in this case.
            \item $a_7\leq 0.27$. \\
            File \texttt{0CD} can examine the first twelve terms $a_1,\dotsc,a_{12}$ to derive a contradiction.
        \end{enumerate}
        \item \label{sim:0D} $0.9\leq a_1+a_2+a_3+a_4\leq 1.1$ and $a_1 \leq 0.3$. \\
        As we need to look up to $a_{16}$, this is again a difficult case. Files \texttt{0D} and \texttt{0DA} successively improve our bounds, with the latter able to prove that $|a_i-1/4|\leq 0.03$ for each $1\leq i\leq 16$. 
        This allows us to apply \Cref{claim:case-d-claim-1} and deduce that $a_1+a_2+a_3+a_4< 1$, which is then enforced in file \texttt{0DAA}. This file, along with \texttt{0DAAA}, then continue to narrow down the acceptable ranges, until file \texttt{0DAAAA} can then prove that $|a_i-1/4|\leq 0.00045$ for each $1\leq i\leq 16$.
        \item \label{sim:0E} $a_1+a_2 \geq 0.9$ and $0.6\leq a_1\leq 0.8$. \\
        File \texttt{0E} proves in this case that there are numbers $c_i\in\set{0,1/3,2/3}$ such that $|a_i-c_i|\leq 0.002$ for each $1\leq i \leq 6$. Looking at the conditions of this case, we can deduce that $c_1=2/3$ and $c_2=\dots=c_6=1/3$.
        \item \label{sim:0F} $0.9\leq a_1+a_2+a_3\leq 1.1$ and $0.45 \leq a_1\leq 0.7$. \\
        Files \texttt{0F} and \texttt{0FA} derive bounds on $a_1,\dotsc,a_{13}$ which allow us to prove that there are numbers $c_i\in\set{0, 0.25,0.5}$ such that $|a_i-c_i|\leq 0.03$ for each $1\leq i \leq 13$.
        Deducing that $c_1=0.5$ and $c_2=\dots=c_{13}=0.25$, we may apply \Cref{claim:case-f-claim-1} and \Cref{claim:case-f-claim-2} to enforce that $a_1+a_2+a_3\leq 1$ and $a_2+a_3+a_4+a_5\leq 1$, which is done in file \texttt{0FAA}. 
        This file, along with file \texttt{0FAAA} prove successively more precise bounds, and the latter is able to prove that, for $c_i$ as above, $|a_i-c_i|\leq 0.00035$ for each $1\leq i \leq 13$.
        \item \label{sim:0G} $1.15\leq a_1+a_2+a_3\leq 1.3$ and $0.45\leq a_2$ and $a_3\leq 0.3$. \\
        File \texttt{0G} proves that for some numbers $c_i\in\set{0, 0.25, 0.5}$, we have $|a_i-c_i|\leq 0.05$ for each $1\leq i\leq 10$. 
        Deducing that $c_1=c_2=0.5$, and $c_3=\dots=c_10=0.25$, we may apply \Cref{claim:case-g-claim-1} and \Cref{claim:case-g-claim-2} to find that $a_2+a_3+a_4\leq 1$ and $a_3+a_4+a_5+a_6\leq 1$.
        These bounds are then enforced in file \texttt{0GA}, and then file \texttt{0GAA} iterates on the resulting bounds to deduce that, for $c_i$ as above, $|a_i-c_i|\leq 0.00045$ for each $1\leq i \leq 10$.
        \item \label{sim:0H} $0.9\leq a_1+a_2+a_3+a_4\leq 1.1$ and $0.3\leq a_1\leq 0.8$. \\
        File \texttt{0H} derives a contradiction. Note that case H of \Cref{sec:lowther-conjecture} is considered a part of \cref{sim:0B} in this section.
        \item \label{sim:0I} $1.1\leq a_1+a_2+a_3\leq 1.2$ and $a_1\leq 0.4$. \\
        File \texttt{0I} splits this case into two further subcases, and then files \texttt{0IA}, \texttt{0IAA}, and \texttt{0IB} derive contradictions in both subcases.
        \item \label{sim:0J} $0.9\leq a_1 + a_2$ and $0.51\leq a_1\leq 0.6$ and $a_2\leq 0.45$. \\
        File \texttt{0J} derives a contradiction.
        \item \label{sim:0K} $1.2\leq a_1+a_2+a_3\leq 1.3$ and $a_2 \leq 0.4$ and $a_3\geq 0.35$. \\
        Files \texttt{0K} and \texttt{0KA} derive a contradiction.
        \item \label{sim:0L} Everything not covered by one of the above cases. \\
        Files \texttt{0L}, \texttt{0LA}, and \texttt{0LAA} derive a contradiction.
    \end{enumerate}
\end{enumerate}

\end{document}